\documentclass{amsart}

\def\tableskip{\medskip}

\theoremstyle{definition}
\newtheorem{algorithm}{Algorithm}
\theoremstyle{remark}
\newtheorem{problem}[algorithm]{Problem}
\theoremstyle{plain}
\newtheorem{lemma}[algorithm]{Lemma}

\usepackage{enumerate}
\usepackage{multirow}
\usepackage[OT4]{fontenc}

\author{Witold Jarnicki}
\address{Institute of Mathematics, Jagiellonian University, Reymonta 4, 30--150 Krak\'ow, Poland}
\email{Witold.Jarnicki@im.uj.edu.pl}
\curraddr{Google R\&E, Krak\'ow, Krupnicza 16, 31--123 Krak\'ow}
\email{witoldjarnicki@google.com}

\author{Maciej \.Zenczykowski}
\address{Theoretical Computer Science Department, Jagiellonian University, Gronostajowa 3, 30--387 Krak\'ow, Poland}
\email{maze@tcs.uj.edu.pl}
\curraddr{Google, 1600 Amphitheatre Parkway, Mountain View, CA 94043}
\email{maze@google.com}

\title{On a property of the number 977731833235239280}
\thanks{%
Part of the computation was done on High Performance Computers of The Academic Computer Centre CYFRONET AGH:
SGI Altix 3700 (baribal, grant number \texttt{MNiSW/SGI3700/UJ/097/2007}),
IBM BladeCenter HS21 (mars, grant number \texttt{MNiSW/IBM\_ BC\_ HS21/UJ/097/2007}),
and
PC RackSaver cluster (zeus, grant number \texttt{MNiSW/PC\_ RS/UJ/097/2007}).
}
\thanks{%
The authors obtained a significant increase of computing speed by using their program to test newly installed Google machines. 
}

\begin{document}

\begin{abstract}
We solve a theoretical arithmetics problem stated by Wac{\l}aw Sierpi\'nski.
The problem has remained open for a couple of decades.
\end{abstract}

\maketitle

\section{Introduction}

Consider the following puzzle from number theory, presented almost 50 years ago.

\begin{problem}[W.~Sierpi\'nski]\label{one}
Find a composite number such that it remains composite after altering any two digits in its
decimal representation.
\end{problem}

We found the problem in \cite{m1}. It was stated as a puzzle for the readers.
It turned out (see \cite{m2}) that no one had solved it. 

A ternary version of the problem was stated in \cite{m2} with a solution given in
\cite{m3}. The basic idea was that altering any two digits in $40=1111_3$ keeps it even and different from $2$. Observe that
such an assumption (changing exactly two digits) makes the problem easy to solve (see Table \ref{exactly_two}). Therefore, it 
is reasonable to assume the following statement.

\begin{problem}[W.~Sierpi\'nski]\label{two}
Find a composite number such that it remains composite after altering at most two digits in its
base $b$ representation.
\end{problem}
\begin{table}[ht!]
\begin{tabular}{|r||r|r||r|r|}
\hline
\multirow{2}*{base}&\multicolumn{2}{|c||}{allowed}&\multicolumn{2}{|c|}{not allowed}\\
\cline{2-5}
&solution&decimal&solution&decimal\\
\hline
\hline
$ 2$&1010100&    84&    1001&9\\
\hline
$ 3$&   1111&    40&      11&4\\
\hline
$ 4$&20130000&34560&12321230&28268\\
\hline
$>4$&       4&    4&       4&    4\\
\hline
\end{tabular}

\tableskip

\caption{Minimal solutions to Problem \ref{one}, 
depending on whether one is allowed to change the most significant digit to zero or not.}
\label{exactly_two}
\end{table}
Another question one may ask is whether it is permitted to change the most significant digit to zero. However,
disallowing that gives little help in finding the number (see Tables \ref{zero} and \ref{no_zero}), so it is safer to assume it may 
be done. 

The problem was
investigated from the theoretical side in \cite{s}. There the author shows, that there are infinitely many solutions to
Problem \ref{two} (for any $b$), provided that Erd\H os's ``favorite'' conjecture on covering 
systems of congruences (see \cite{e}) is true. However, since Erd\H os's conjecture is open, so remained the problem.

\section{Main results}

We present the solution to Problem \ref{two} (bases 2--10) in Tables \ref{zero} and \ref{no_zero}. 
%
\begin{table}[ht!]
\begin{tabular}{|r||r|r|}
\hline
base&solution&decimal\\
\hline
\hline
 2&           1010100&                84\\
\hline
 3&           2200100&              1953\\
\hline
 4&          20130000&             34560\\
\hline
 5&        3243003420&           7000485\\
\hline
 6&       55111253530&         354748446\\
\hline
 7&     5411665056000&       77478704205\\
\hline
 8&    33254100107730&     1878528135128\\
\hline
 9&   210324811482600&    48398467146642\\
\hline
10&977731833235239280&977731833235239280\\
\hline
\end{tabular}

\tableskip

\caption{Solutions to Problem \ref{two} found for bases 2--10 (changing the most significant digit to zero is allowed). 
They are known to be minimal for bases 2--9.}
\label{zero}
\end{table}
\begin{table}[ht!]
\begin{tabular}{|r||r|r|}
\hline
base&solution&decimal\\
\hline
\hline
 2&           1010100&                84\\
\hline
 3&           2200100&              1953\\
\hline
 4&          12321230&             28268\\
\hline
 5&         324322330&           1401590\\
\hline
 6&       43040303150&         273241578\\
\hline
 7&     5411665056000&       77478704205\\
\hline
 8&    33254100107730&     1878528135128\\
\hline
 9&   210324811482600&    48398467146642\\
\hline
10&977731833235239280&977731833235239280\\
\hline
\end{tabular}

\tableskip

\caption{Solutions to Problem \ref{two} found for bases 2--10 (changing the most significant digit to zero is disallowed). 
They are known to be minimal for bases 2--9.}
\label{no_zero}
\end{table}
\section{Motivation and methods}

Since there are many ways to modify a number by altering two of its digits, it initially seems impossible to find a solution 
to the
problem. Therefore, a natural way to approach the puzzle is to consider it in bases smaller
than $10$ --- in such situation the number of ways of changing a number is much smaller. 

We managed to solve the problem for bases between $2$ and $9$ using a small
``grid'' of computers at our university. The computation lasted several weeks.
Looking at the results gives grounds to suppose that solving the problem for base $k+1$
requires about $100$ times the computing power needed to solve it for base $k$.
This estimate inspired us to use grid computing to solve the decimal case.

To present the main idea behind the solution method, consider the following lemma.

\begin{lemma}\label{lem}
Let $n$ contain at least $5$ nonzero digits in base $b$. Assume that $n$ is a solution to Problem \ref{two}. 
Define $\widetilde n:=b\lfloor n/b\rfloor$. Then $\widetilde n$ is a solution to Problem \ref{two}.
\end{lemma}

\begin{proof}
The number $\widetilde n$ is composite, because $b|\widetilde n$ and $1<b<\widetilde n$. For the same reason it will remain
composite whenever we decide to leave the least significant digit intact. To complete the proof, observe that altering two
digits of $\widetilde n$, one of which is the least significant one, leads us to a number that is obtained by exactly the
same alteration of $n$.
\end{proof}

Observe that, for $n$ having $1$ nonzero digit, either its least significant digit is zero (then $\widetilde n=n$ and
Lemma \ref{lem} holds trivially), or nonzero (then $n$ can be transformed into $1$, which is not composite, and Lemma \ref{lem} holds
trivially). 

Additionally, for $n$ having $4$ nonzero digits, the only situation when
Lemma \ref{lem} fails is when the second least significant digit of $n$ is $1$.

This means that, when looking for a solution to Problem \ref{two}, it is sufficient to consider a relatively small 
number of cases where $n$ has $2$, $3$, or $4$  nonzero digits (in the latter 
case with second least significant digit being $1$)
and the numbers divisible by $b$. The former can be done directly.
The latter is accomplished by the following algorithm.

\begin{algorithm}\label{main}
Given integers $b\geq 2$, $e\geq 0$, $u\geq 0$, to find all solutions
$n\in b\mathbb Z\cap[ub^{e+1},(u+1)b^{e+1})$ to Problem \ref{two}:
\begin{enumerate}[1.]
\item\label{sieve_begin} 
For each integer $k\in[0,b^{e+1})$:
\item\label{sieve_end}
\quad If $ub^{e+1}+k$ is prime: set $a_k\longleftarrow1$; otherwise set $a_k\longleftarrow0$.
\item 
For each integer $k\in[0,b^e)$:
\item
\quad If any of $a_{bk},a_{bk+1},\dots,a_{bk+b-1}$ is $1$: set $b_k\longleftarrow1$; otherwise set $b_k\longleftarrow0$.
\item
For each integer $k\in[0,b^e)$:
\item\label{solution_begin}
\quad Set $c_k\longleftarrow1$.
\item\label{cross_begin}
\quad For each integer $\ell\in[0,b^e)$ with $k$ differing from $\ell$ at at most one digit:
\item\label{cross_end}\label{solution_end}
\quad\quad If $b_\ell=1$: set $c_k\longleftarrow0$.
\item
For each integer $k\in[0,b^e)$ with $c_k=1$: 
\item\label{check}
\quad Check directly if $b(ub^e+k)$ is a solution and output if it is.
\end{enumerate}
\end{algorithm}

The steps \ref{sieve_begin}--\ref{sieve_end} require finding all prime numbers in an interval. We have used the following 
algorithms, all with running time close to $O(b^e)$:
\begin{enumerate}[(a)]
\item a straightforward implementation of the sieve of Eratosthenes,
\item the implementation of the sieve of Atkin from \cite{b},
\item our own implementations of the sieve of Atkin with $W=12$, $W=60$, and $W=420$ (see \cite{ab} for details).
\end{enumerate}
The author of \cite{b} claims that the program works for primes up to $10^{15}$. The code is pretty complicated, so we 
could not figure out whether it works past that boundary. That is why we decided to create our own implementation.
Surprisingly, the version $W=12$ worked best for really large numbers.

The steps \ref{cross_begin}--\ref{cross_end} can be done in $O(eb)$ time. Consequently, we implemented the whole loop 
\ref{solution_begin}--\ref{solution_end} in $O(eb^e)$ time.

In practice, the step \ref{check} is involved only for a couple values of $k$, so its influence on running time is negligible.

\bigskip

If we are interested in finding the smallest solution, it is enough to call Algorithm \ref{main} for fixed $b,e$ and 
sequential units $u$. However, due to limited computational resources, it is more important to find any solution rather than 
to prove that it is minimal. Therefore, we first scanned the units, for which the probability of finding a solution is high.

The method we used is neither strict nor formal, but worked in practice. For the sake of the estimation, we assumed that 
primality of numbers is a result of a sequence of independent random experiments. The probablity of a number in unit $u$ being
prime is $(\pi((u+1)b^{e+1})-\pi(ub^{e+1}))b^{-(e+1)}$. Using the approximation $\pi(n)\approx n/\ln(n)$ we obtained the 
approximate value of the
probability $p(u)$ of a number from block $u$ being a solution to Problem \ref{two}. Then, we considered the blocks in 
order of decreasing $p(u)$. The computation done is summed up in Table \ref{computation}.
\begin{table}[ht!]

\begin{tabular}{|c|c||c|c|}
\hline
$b$&$e$&units scanned&solution unit\\
\hline
\hline
2&6&0&0\\
\hline
3&6&0&0\\
\hline
4&7&0&0\\
\hline
5&9&0&0\\
\hline
6&10&0&0\\
\hline
7&10&0--39&39\\
\hline
8&10&0--218&218\\
\hline
9&10&0--1542&1542\\
\hline
10&9&0--97773183&97773183\\
\hline
\end{tabular}

\tableskip

\caption{Summary of units considered}
\label{computation}
\end{table}

\end{document}